\documentclass[12pt,oneside,mathscr]{amsart}
\usepackage{amsthm,amsmath,amssymb,euscript,amscd,a4wide}
\usepackage{bigstrut}
\usepackage[all]{xy}
\pdfoutput=1

\newtheorem{thm}{Theorem}[section]
\newtheorem{cor}[thm]{Corollary}
\newtheorem{lemma}[thm]{Lemma}
\newtheorem{prop}[thm]{Proposition}

\theoremstyle{definition}
\newtheorem{definition}[thm]{Definition}
\newtheorem{note}[thm]{Notes}

\theoremstyle{remark}
\newtheorem{remark}[thm]{Remark}
\newtheorem{notation}[thm]{Notation}

\numberwithin{equation}{section}

\newcommand{\thmref}[1]{Theorem~\ref{#1}}

\newcommand{\lemref}[1]{Lemma~\ref{#1}}
\newcommand{\propref}[1]{Proposition~\ref{#1}}
\newcommand{\seqref}[1]{sequence~\ref{#1}}
\newcommand{\corref}[1]{Corollary~\ref{#1}}

\newcommand{\C}{{\mathbb C}}

\newcommand{\PP}{{\mathbb P}}
\newcommand{\OO}{{\mathscr O}}
\let\cross=\times
\let\tensor=\otimes

\newcommand{\spec}{\operatorname{spec}}
\newcommand{\Jac}{\operatorname{Jac}}
\newcommand{\irf}[1]{\hat\Phi^{#1}}
\newcommand{\rf}[1]{\Phi^{#1}}
\newcommand{\muk}{\operatorname{\Phi}}
\newcommand{\dmuk}{\operatorname{\hat\Phi}}
\newcommand{\Hilb}{\operatorname{Hilb}}

\newcommand{\Pic}{\operatorname{Pic}} 
\newcommand{\Ext}{\operatorname{Ext}}

\newcommand{\coker}{\operatorname{coker}}

\newcommand{\ch}{\operatorname{ch}}
\newcommand{\rk}{\operatorname{rk}}

\newcommand{\T}{{\mathbb T}}
\renewcommand{\P}{{\mathscr P}}

\newcommand{\dT}{{\hat{\mathbb T}}}
\newcommand{\flip}{(-1_{\T})^*}

\newcommand{\invmuk}{\operatorname{\hat\Phi}}

\newcommand{\xhat}{{\Hat x}}
\newcommand{\yhat}{{\Hat y}}
\newcommand{\zhat}{{\Hat z}}

\newcommand{\dhat}{{\Hat d}}

\newcommand{\Lhat}{{\Hat L}}
\newcommand{\Lshat}{\widehat{L^2}}
\newcommand{\I}{{\mathscr I}}

\newcommand{\compo}{\raise2pt\hbox{$\scriptscriptstyle\circ$}}
\renewcommand{\leq}{\leqslant}
\renewcommand{\geq}{\geqslant}
\newcommand{\lra}{\longrightarrow}

\newcommand{\lRa}[1]{\>{\buildrel {#1}\over\longrightarrow}\>}

\newcommand{\sing}{\operatorname{sing}}

\newcommand{\pihat}{\Hat\pi}
\newcommand{\supp}{\operatorname{supp}}
\newcommand{\length}{\operatorname{length}}

\begin{document}

\title[]{A Fourier-Mukai Approach to the Enumerative Geometry of Principally Polarized Abelian Surfaces}
\author[]{Antony Maciocia}
\address{Department of Mathematics and Statistics\\
The University of Edinburgh\\
The King's Buildings\\ Mayfield Road\\ Edinburgh, EH9 3JZ.\\}
\email{a.maciocia@.ed.ac.uk}
\thanks{} 
\date{\today}

\begin{abstract}
We study twisted ideal sheaves of small length on an irreducible principally
polarized abelian surface $(\T,\ell)$. Using Fourier-Mukai techniques
we associate certain jumping schemes to such sheaves and completely
classify such loci. We
give examples of applications to the enumerative geometry of $\T$ and
show that no smooth genus $5$ curve on such a surface can contain a
$g^1_3$. We also describe explicitly the singular divisors in the 
linear system $|2\ell|$.
\end{abstract}

\maketitle

\section*{Introduction}{}
It is an old problem of algebraic geometry to describe the family of
curves in a given linear system which go through a certain number of points.
For projective spaces, this problem is essentially solved. The answer is given
by the Pl\"ucker formulae and the reciprocity formulae (see \cite[\S2.4
and \S5.4]{GriffHarris}). Our aim, in this
paper, is to extend some of these results to a principally polarized Abelian
surface $(\T,\ell)$ over the complex numbers.
We shall find that, for a generic set of points, the family of curves has the
dimension we would expect from a na\"\i ve dimension count. We aim to solve
the problem of determining these families explicitly for 
divisors in translates of $|2\ell|$, this includes the non-generic
cases. This will entail studying the case $|\ell|$ as well. The unique
divisor $D$ in $|\ell|$ and its translates will be viewed as ``lines''
on the abelian surface. There is a duality for these divisors in that any two intersect
in two points (with multiplicity) while any two distinct points lie on exactly
two of these ``lines''. We use the term ``collinear'' to mean that a
0-scheme lies on a ``line''.

Our main tools are sheaf theory and the Fourier-Mukai transform which provides
exactly enough information to determine the incidence of points on
divisors. We proceed in an inductive way: we relate the properties of
the twisted ideal sheaf
$L^2\tensor\I_X$ to those of $L^2\tensor \I_{X'}$ for $X'\subset
X$. We then give an exhaustive treatment of $L^2\tensor\I_X$ for each
small value of $|X|$ in turn. 

The paper is organized as follows. Section one sets up some convenient
notation for zero-dimensional subschemes. The second section describes some of
the basic properties of the linear systems $|\ell|$ and $|2\ell|$
including a description of the reducible divisors in $|2\ell|$. These
are just the $\Theta$ and 2$\Theta$ linear systems and the reducible
divisors are very well understood; forming the cornerstone of the
theory of principally polarized abelian surfaces. But we recall some
of the key facts translated into the language of sheaves. Section
three is a brief overview of the Fourier-Mukai transform. In section four we
shall define the basic objects of study: the cohomology jumping schemes
associated to a zero-dimensional subscheme of $\T$. These
describe the translates of $|2\ell|$ which contain extra divisors which
go through $X$. 

In sections five to nine we study the cases of length 1, 2, 3, 4 and 5 subschemes in
detail. Each of these requires special treatment and a full analysis of each
is required before proceeding to the next. In section 10 we can treat the
general case. The conclusions are summarized in the two tables of the appendix
to the paper. In section 11 we look briefly at stability questions for
the Fourier transforms of the twisted ideal sheaves and in section 12
we look at two applications, one to answer the question of whether
smooth divisors in $|2\ell|$ have a $g^1_3$ or not (they don not, as
we shall see) and to answer a classical question to compute the
locus of singular divisors in $|2\ell|$. The idea is to use the
information about how non-reduced 0-schemes lie on divisors to
detect singularities in the divisors of $|2\ell|$.

Our motivation for studying these questions come from the need to compute
moduli spaces of stable and semistable sheaves over these tori. These moduli
spaces provide an endless source of (relatively) easy to compute
hyper-K\"ahler manifolds which lie in a good deformation family parametrized
by the moduli space of flat hyper-K\"ahler tori. Some of these moduli spaces
can be related to the space of 0-dimensional subschemes of the torus and the
stability and local-freeness properties of the associated sheaves is
determined by the incidence of the 0-subscheme on divisors from $|L^n|$. 
As an example of this, one can apply the results of sections eight, nine and ten
to determine the
moduli space of stable sheaves with Chern characters $(2,0,-2)$ and $(2,0,-3)$. It
turns out that the standard algebraic compactification of the latter moduli
space is isomorphic to $\Hilb^6\T\cross\T$. The former is isomorphic to
another compactification of the bundle of Jacobians over the space of
effective divisors in all translates of $|2\ell|$ (see
\cite{Mac0}). It also provides us with a new irreducible hyperK\"ahler
manifold (see \cite{OGradyNewSymp}). 
This paper forms the basis for a research programme which aims to give
complete descriptions of a variety of moduli spaces of sheaves and
more generally Bridgeland stable objects on a principally polarized
abelian surface. The importance of moduli spaces arising from Hilbert
schemes of points can be seen in \cite{Yosh01} where it is shown that
each fine moduli spaces of stable sheaves are birational to
some $\Hilb^n\T\times\dT$. 

The results we find in this paper are also useful to study moduli of Bridgeland stable
objects with the same Chern character and in a subsequent work we
compute the explicit wall crossing behaviour of such moduli spaces. In
this paper, we operate at a more elementary level and avoid the use of
Bridgeland stability but note that in section 11 a fuller treatment is
best completed using the derived category. This has been done in
\cite{MacMeachan} where the Bridgeland stable moduli spaces for Chern
character $(1,2\ell,n)$ are computed. It is also possible to give an
account of $n$ very ampleness for powers of $L$ using these
techniques. This will be the subject of a future article.

\section{The Hilbert scheme of 0-subschemes}{}
\label{s:hilbscheme}
Let $S$ be a smooth complex surface and let $\Hilb^nS$ denote the Hilbert
scheme of length $n$ 0-dimensional subschemes of $S$.

It is useful to introduce the following notation.
A general 0-dimensional scheme will be denoted by $X$. Its length will be
denoted $|X|$.
\begin{alignat*}2
P&=\{p\}&&\text{ a single point,}\\
Q&=\{p,q\}&&\text{ a length 2 0-scheme,}\\
&&&\text{\ \ if $p=q$ then write $\{p,t\}$, where $t\in\PP T_pS$}\\
Y&=\{p,q,y\}&&\text{ a length 3 0-scheme,}\\
Z&=\{p,q,y,z\}&&\text{ a length 4 0-scheme,}\\
W&=\{p,q,y,z,w\}&&\text{ a length 5 0-scheme}
\end{alignat*}

The ideal sheaf associated to $X\subset S$ will be denoted $\I_X=\ker(\OO_S\to
\OO_X)$. The notation $P,Q,Y,Z,W$ is used to avoid double subscripts. It
also reflects the fact that the properties of 0-schemes of different small
lengths vary widely with respect to a given linear system whereas large length
0-schemes all behave similarly. As a mnemonic for remembering which
are which note that the drawn letters have their length number of ``points''
on them (at the ends of lines or at acute angles). 

Suppose that $D\subset S$ is a curve on a variety then we
shall be concerned with questions of the following type: given a 0-dimensional
subscheme $X\subset S$, which curves in a given linear system or systems
containing $D$ contain $X$? In other words, we want to understand
$H^0(\OO(D)\tensor\I_X)$.

Suppose $X$ is non-reduced and
contains precisely one closed point so that it is of the form $\spec A$ for
some Artin local ring $A$. These will take the form $A=\C[\epsilon,\eta]/I$
for some ideal $I$.
\begin{notation}\label{n:Yx}
We will have need to distinguish the three possibilities
for the isomorphism class of $Y$ (as a scheme)
when $\supp(Y)=\{p\}$. We shall denote these as follows:
\begin{equation*}
Y_c=\spec\frac{\C[\epsilon]}{(\epsilon^3)},\qquad
Y_d=\spec\frac{\C[\epsilon,\eta]}{(\epsilon-\eta^2,\epsilon\eta)}
\quad\text{and}\quad
Y_e=\spec\frac{\C[\epsilon,\eta]}{(\epsilon^2,\eta^2,\epsilon\eta)}
\end{equation*}
Note that the ideal giving rise to the first of these is regular while the
ideals of the other two are not. This means, in particular, that the dimension
of the fibre of $\I_Y$ at its support is 2 for $Y_c$ and 3 for $Y_d$ and $Y_e$.
\end{notation}
Recall that a node is a transverse
intersection, a cusp is a locally irreducible double point of the curve
and a tacnode is a
double point which is not transverse. In terms of the classification of
surface singularities given in \cite[\S II.8]{BPVdV}: a node has type $A_1$,
cusps have types $A_{2n}$ and tacnodes have types $A_{2n+1}$. By a
\textit{simple} tacnode we mean an $A_3$ singularity. 

\section{The linear systems $|\ell|$ and $|2\ell|$}{}
We shall now suppose that $(\T,\ell)$ is a principally polarized Abelian surface
and $(\dT,\hat\ell)$ is the dual torus. We let  $\T_2$ denote the set
of points of order 2 in $\T$. We choose a symmetric line 
bundle $L$ in the class $\ell$. From now on, we will use the
notation $L^n$ to denote the
$n$th tensor product of $L$ with itself (rather than $L^{\tensor n}$).
Note that $h^0(L)=1$ and we write $D_L$
for the zero set of the unique non-zero holomorphic section of $L$.
In this and subsequent sections we adopt the following
convention: $D_u$ will denote the translation of $D_L$ by $u\in\T$. This makes
the notation more concise and also allows us to write $T_r D_u$ to mean the
tangent space at $r$ of the translate of $D_L$ by $u$. Observe the irritating
anomaly that $D_u\in |\tau_{-u}^*L|$, where $\tau_x:\T\to\T$ denotes
translation by $x$.

 If $D_L$ is irreducible then
$\T=\Jac(D_L)$ and $L$ corresponds to the $\theta$ polarization and then $D_L$
must be smooth (see e.g. \cite[Cor 11.8.2]{BirkLange}).

We shall introduce the following terminology.
\begin{definition}
If we have  $X\in\Hilb^n\T$ and $X\subset D_u$ for some $u\in\T$ then we say that $X$ is \emph{collinear}.
\end{definition}
This slightly unorthodox use of `collinear' makes some sense because the
divisors $D_u$ play a similar role to lines on $\C P^2$. In fact, we
shall see that any two ``lines'' intersect in exactly two points (up
to multiplicity) and dually, that any two distinct points are contained in
exactly two ``lines''.

Turning now to $|2\ell|$, observe that $h^0(L^2)=4$ and recall that $L^2$ is
base-point free and so, by Bertini, the generic element $D\in|2\ell|$ is smooth.
The reducible divisors in $|2\ell|$ are given by the following.
\begin{lemma}\label{l:DinL2}
If $D\in|2\ell|$ is reducible then 
 $D=D_x+D_{-x}$.
\end{lemma}
The proof is a straightforward exercise (see \cite[Chapter 10]{BirkLange}).

We use the notation $Ks(\T)=\T/\pm$ to
denote the singular Kummer variety. This has sixteen
singular points and they all have type $A_2$ (in terms of the classification
of surface singularities given in \cite[III.3]{BPVdV}). 
%By the Theorem of the Square (see \cite[p59]{Mum})
There is a canonical family of divisors in
$|\Lhat^2\tensor\P_\gamma|$ given by $D_u+D_y$, where
$u+y=\gamma$. This family is parametrized by $Ks_\gamma(\dT)$. We shall call
the divisors \textit{Kummer divisors}.
The map $Ks_\gamma(\dT)\to|\Lhat^2\tensor\P_\gamma|$ is
given by $[\alpha]\mapsto D_\alpha+D_{\gamma-\alpha}$, where
$\alpha\sim\alpha'$ if and only if $\alpha=\gamma-\alpha$, so that
$Ks_\gamma(\dT)=\dT/{\sim}{}\cong Ks(\dT)$ isomorphic to
the image of the above map. 
\begin{notation}
We denote the Chern characters of sheaves on $\T$ by $(r,c_1,\chi)$,
where $\chi(E)=\frac{1}{2}c_1(E)^2-c_2(E)$ and $r=\rk(E)$.
\end{notation}
The following lemma (which works over any smooth surface) will be used several times.
We let $\mathcal{L}(D)$ denote the line bundle associated to a divisor.
\begin{lemma}\label{l:SubSch}
Suppose $X\subset D$ is a 0-subscheme of an effective divisor $D$ on a surface
$S$. Let $R$ be a line bundle over $S$ and $A=R\I_X/R\mathcal{L}(D)^*$.
Then $A$ contains no subscheme supported in codimension 2.
\end{lemma}
\begin{proof}
If $T\subset A$ is supported in codimension 2 then
$\ch_2(A/T)\leq\ch_2(A)$ and so $\ch_2(K)\geq\ch_2(R\mathcal{L}(D)^*)$, where
$K=\ker(R\I_X\to  A/T)$. But $K=R\mathcal{L}(D)^*\I_{X'}$ for some 0-scheme $X'$ and so
$\ch_2(K)=\ch_2(R\mathcal{L}(D)^*)-|X'|$. Hence $|X'|=0$ and $T=0$.
\end{proof}

Finally, we will need to describe torsion-free sheaves of rank 1 over
reducible curves $D_1+D_2\subset \T$. If $T$ is such a sheaf then we can
consider its restrictions $T_i=T|_{D_i}$. Define the \textit{degree} of $T$ via
the Riemann-Roch formula for embedded curves (see \cite[II.3]{BPVdV}) to be
\[\deg(T)=\chi(T)-\chi(\mathcal{L}(D_1+D_2)).\] 
Define the \textit{restriction type} or \textit{degree} of $T$ to be
$(\deg(T_1),\deg(T_2))$. If the restriction type of $T$ is $(n_1,n_2)$  and
the singularity scheme $X$ of $T$ has length $r$ then it is easy to see that
$n_1+n_2=\deg(T)-r$. Note that $X\subset D_1\cap D_2$ and so $r\leq |D_1\cap
D_2|$ if $D_1$ and $D_2$ have no common components. This generalizes in the 
obvious way to more than two components and to multiple curves.

\section{The Fourier-Mukai Transform}

One of the most useful tools in studying
sheaves and bundles over tori is the Fourier-Mukai transform (see e.g.
\cite{MukDual} or \cite{MukFF} for the original treatment and
\cite{HuyBook} or \cite{BBHBook} for more up to date treatments).
This takes the form of a functor $\muk:D(\T)\to D(\dT)$ between the derived
categories of complexes of coherent sheaves. It is defined by
$E\mapsto\mathbf{R}\hat\pi_*(\pi^*E\tensor\P)$, where $\P$ denotes the Poincar\'e
line bundle over $\T\cross\dT$ and $\pi$ and $\pihat$ are the projection maps
to $\T$ and $\dT$ respectively. We denote the
fibres of $\P$ over $\T\cross\{\xhat\}$ by $\P_\xhat$ which gives rise to the
isomorphism $\dT\cong\Pic^0\T$. A particularly important result concerning the
Mukai transform is the fact that if $\rf j(E)=0$ for all $j>m$ then the fibres
of $\rf m(E)$ are given (canonically) by $H^m(E\tensor\P_\xhat)$.
\begin{definition}\label{d:WITIT}
Following
Mukai, we say that a sheaf $E$ over $\T$ satisfies WIT$_i$ if $\rf j(E)=0$ for
all $j\neq i$ and write $\hat E$ for $\rf iE$. We also say that $E$ satisfies
IT$_i$ if $H^j(E\tensor\P_\xhat)=0$ for all $\xhat\in\dT$ and $j\neq i$.
\end{definition}
Note that ample line bundles satisfy IT$_0$. There is an unavoidable conflict
of notation with $\hat L$. This will always denote the dual polarization which
equals $(\rf0L)^{-1}$. Notice that the fibres of the projective bundle $\PP\Lshat$
are given canonically by the linear systems $|L^2\tensor\P_\xhat|$ as $\xhat$
varies over $\dT$.

The inverse of $\Phi$ is $(-1)^*\hat\Phi[2]$ where
$\dmuk(E)=\mathbf{R}\pi_*(\hat\pi^* E\tensor \P)$. 
For practical purposes the fact that $(-1)^*\hat\Phi[2]$ is the
quasi-inverse of $\muk$ can be viewed as the
following first-quadrant spectral sequence
\[
E^{p,q}_2\>\Rightarrow\>\left\{
\begin{array}{ll}\flip E&\qquad\text{ for $p+q=2$},\\
0&\qquad\text{ otherwise,}\end{array}\right.
\]
with $E^{p,q}_2=\dmuk^p(\muk^q E)$. 
But $\Phi^2(E)$ is IT$_0$ while $\Phi^0(E)$ is WIT$_2$.
Hence, the entire information content of this spectral sequence
is contained in the exact sequences
\[0\lra\mathcal{D}\lra(-1_{\T})^*E\lra \dmuk^0(\muk^2(E))\lRa{d_2}\dmuk^2(\muk^1(E))\lra0\]
and
\[0\lra \dmuk^0(\muk^1(E))\lRa{d_2}\dmuk^2\muk^0(E))\lra\mathcal{D}\lra \dmuk^1(\muk^1(E))\lra0,\]
where $\mathcal{D}$ is an unknown. Moreover, if any one of
$\dmuk^2(\muk^0(E))$, $\dmuk^1(\muk^1(E))$, $\dmuk^0(\muk^2(E))$ is zero then 
we can eliminate $\mathcal{D}$.

We list some useful properties of $\muk$.
\begin{prop}\label{p:Muklem}
(see \cite{MukDual})
\begin{enumerate}
\item $\OO_X$ satisfies IT$_0$ and $\rf0(\OO_X)=H_X$, a homogeneous
bundle\\
and $\rf0(\OO_x)=\P_x$.
\item $\P_\xhat$ satisfies WIT$_2$ and $\rf2(\P_\xhat)=\OO_{-\xhat}$
\item If $E$ satisfies WIT then so does $\tau^*_xE$ with transform $\hat
E\tensor\P_{-x}$.
\item If $\ch(E)=(r,c,\chi)$ then $\ch(\muk(E))=(\chi,-\hat c,r)$.
\end{enumerate}
\end{prop}
\begin{remark}
Observe that $\PP(\widehat{L^i})$ is
flat as a projective bundle. Hence, $\widehat{L^i}$ admits an irreducible
projectively flat connection and so $\widehat{L^i}$ is $\mu$-stable. This is
proved in a different way by Kempf (see \cite{kempf}{Thm 3}).
By $\mu$-stable we mean the stability of Mumford-Takemoto: $E$ is
{\em $\mu$-stable}
if $E$ is torsion-free and for all subsheaves $F\subset E$ with
$E/F$ torsion-free we have $d(F)/r(F)<d(E)/r(E)$. We obtain
{\em $\mu$-semistability} by replacing $<$ by $\leq$. Homogeneous bundles can
be characterized as $\mu$-semistable sheaves with $c_1=0$ and $c_2=0$.
\end{remark}

\begin{notation}
As a useful shorthand, we shall drop the tensor product sign when no
confusion will arise. We shall also write
$L_\xhat=L\P_\xhat=L\tensor\P_\xhat$.
\end{notation}
We will be interested in $L^i\I_X$ for $i>0$. Since $\muk$ is right
exact we can
apply it to short exact sequences to obtain a long exact sequence. For
example,
\[\muk(0\lra L^i\I_X\lra L^i\lra\OO_X\lra0)\]
gives rise to
\begin{equation}\label{e:RFofstr}
0\to \rf0(L^i\I_X)\to \widehat{L^i}\to
H_X\to\rf1(L^i\I_X)\to0.
\end{equation}
From this sequence we can immediately deduce:
\begin{prop}\label{p:R0FLI}
For all $i>0$ and 0-schemes $X$, $\rf2(L^i\I_X)=0$ and hence
$H^2(L^i\I_X\P_\xhat)=0$ for all $\xhat$. 
\end{prop}
In particular,
$\rf0(L^i\I_X)$ is locally-free and, since $\widehat{L^i}$ is stable we
have $\mu(\rf0(L^i\I_X))<-1$.
Observe that $\chi(L^i\I_X)=i^2-|X|$.

\begin{definition}\label{d:Rij}
Let $R^j_i(X)=\rf j(L^i\I_X)$ for $i>0$ and $j=0,1$.
\end{definition}
We can apply the Mukai spectral sequence to $L^i\I_X$ to obtain a long
exact sequence
\begin{equation}\label{e:MSSLI}
0\to\irf0(R^1_i(X))\to\irf2(R^0_i(X))\to L^i\I_{-X}\to
\irf1(R^1_i(X))\to0.
\end{equation}
Observe that it is impossible for the middle map to be zero unless
$L^i\I_X$ satisfies WIT$_1$.

If $X'\subset X$ with $X\setminus X'=X''$ then the sequence
\[0\lra L^i\I_X\lra L^i\I_{X'}\lra\OO_{X''}\lra0\]
gives rise to
\begin{equation}\label{e:RFofPstr}
0\to R^0_i(X)\to R^0_i(X')\to H_{X''}\to R^1_i(X)\to R^1_i(X')\to0.
\end{equation}
This is particularly useful if $|X''|=1$.

\section{Cohomology Jumping Schemes}{}
Let $p:U\to S$ be a flat morphism of projective varieties. Let $F$ be a sheaf
on $U$ and let $F_r$ denote its fibre over $r\in S$.
Suppose that, for some sheaf $F$,
$R^ip_*F=0$ for $i>1$. Then the fibres of $R^1p_*F$ over $r\in S$ are
canonically isomorphic to $H^1(F_r)$ since $H^i(F_r)=0$ for all $i>1$. 
We want to consider the points of $S$ where the dimension of $H^1(F_r)$ jumps
(up, by semicontinuity). We will consider the situation where $F=
p_1^*(L^i\I_X)\tensor\P$ over $\T\cross\dT$ and $S=\dT$. Since $\chi(F_r)$
does not depend on $r$ we see that $\dim(H^0(L^i\I_X\P_\xhat))
=\dim(H^1(L^i\I_X\P_\xhat))=\dim(\rf1(L^i\I_X)\tensor\OO_\xhat)$. From
\seqref{e:RFofstr} we know that
\[\rf1(L^i\I_X)=\coker(\rf0(L^i\to\OO_X):\widehat{L^i}\to
H_X).\] We can therefore make the following definition.

\begin{definition}
The \textit{cohomology jumping scheme} $S_i(X)$ associated to $L^i\I_X$ is
defined to be the determinantal locus of $\rf0(L^i\to\OO_X)$.
\end{definition}
More often than not, we shall only be interested in the support of $S_i(X)$.
\begin{definition}
Define $\Phi_i(X)\subset \PP\widehat{L^i}$ to be $\{D\mid X\subset D\}$.
Note that $\PP\widehat{L^i}$ is a projective bundle over $\T$ and we denote
the intersection of $\Phi(X)$ with the fibre over $\xhat$ by $\Phi_i(X)_\xhat$.
\end{definition}
Observe that
\[\Phi_i(X)_\xhat\cong \PP H^0(L^i\P_\xhat\I_X).\]
By semicontinuity of cohomology, for generic $\xhat$, $\Phi_i(X)_\xhat\cong
\C P^r$ for some $r$. Then, for such $\xhat$,
\[h^1(L^i\I_X\P_\xhat)=|X|-i^2+r.\]
In
any case, $r\leq\dim\Phi_i(X)_\xhat\leq i^2$. It also follows that the support
of $S_i(X)$ is just $\{\xhat\in\dT\mid \dim\Phi_i(X)_\xhat\geq r+1\}$.

Our aim will be to compute $\Phi_1(X)$ and $\Phi_2(X)$ for any $X$. We do this
by considering separately the cases $|X|=1,2,3,4$ and then extrapolating to
the general case.

\section{$|X|=1$}{}

Consider $S_1(P)$ first. From \seqref{e:RFofstr} we see that $L\I_P$
satisfies WIT$_1$ and its transform is given by $\P_p/\hat L^{-1}=
\OO_{D_{-p}}\P_p$. Hence, $S_1(P)=D_{-p}$. In fact, any degree 0 line
bundle over $D_u$ satisfies WIT$_1$ with transform $L_\beta\I_{-u}$. 

Turning now to $S_2(P)$, observe that $\chi(L^2\I_P)=3$ and hence $\dim
H^0(L^2\I_P)\geq3$. Thus the fibres
of $\Phi_2(P)$ contain $\C P^2$. But $|L^2\P_\xhat|$ are all base-point
free and hence the fibres of $\Phi_2(P)$ must all equal $\C P^2$. In other
words, $S_2(P)=\varnothing$ and $L^2\I_P$ satisfies IT$_0$.

\begin{prop}
For all $P$, $R^0_2(P)=\Phi^0(L^2\I_p)$ are $\mu$-stable.
\end{prop}
\begin{proof}
This follows because $\Phi^0(L^2)$ is $\mu$-stable and $R^0_1(P)$ is a
sub-bundle of $\Phi^0(L^2)$ of slope $-4/3$. Then there are no
integers $a$ and $b$ such that $-4/3\leq 2a/b<-1$ with $0<b<3$.
\end{proof}

\section{$|X|= 2$}{}

Part of the following is well known (see \cite{BirkLange}).
\begin{prop}\label{p:S1Q}
For all $Q\in\Hilb^2\T$, $L\I_Q$ satisfies WIT$_1$ and $R^1_1(Q)=\hat
L_{p+q}\I_{S_1(Q)}$, where $S_1(Q)\in\Hilb^2\dT$. 
\end{prop}
\begin{proof}
Use sequence \ref{e:RFofstr} with $i=1$ and $X=Q$. Then $R^0_1(Q)=0$ or the
sequence
splits. But $H_Q$ satisfies WIT$_2$ whereas $\irf2(R^1_1(Q))=0$ from the Mukai
spectral sequence. Hence, $L\I_Q$ satisfies WIT$_1$.

Now sequence \ref{e:RFofPstr} with $\{p\}\subset Q$ gives
\begin{equation}\label{e:R11Q}
0\lra\P_q\lra R^1_1(Q)\lra\OO_{D_{-p}}\P_p\lra0.
\end{equation}
This cannot split as $\irf2(R^1_1(Q))=0$. Hence,
$R^1_1(Q)$ is torsion-free and so takes the form $\hat L_x\I_{Q'}$ with
$Q'\in\Hilb^2\dT$. From \ref{e:R11Q}, we see that $x=p+q$.
\end{proof}

As an easy corollary of this and sequence \ref{e:RFofPstr}
to give the answer for $S_1(X)$ for any $X$.
\begin{cor}\label{c:genS1}
If $X$ is a 0-scheme and length at least 2 then $L\I_X$ satisfies WIT$_1$
and $R^1_1(X)$ is torsion-free.
\end{cor}
We can view \propref{p:S1Q} more explicitly as follows.
Using the fact that $S_1(P)=D_{-p}$ we see that
$S_1(Q)=D_{-p}\cap D_{-q}$. If $Q$ is reduced we
can deduce that $S_1(Q)=\{l-p,-l'-p\}$, where $p-q=l-l'$ for $l,l'\in D_L$.
Furthermore, $S_1(Q)$
is reduced if and only if $p-q\not\in 2D_L$.

Suppose $Q$ is not reduced and given by
$t\in\PP T_p\T$. We have a degree 2 map
$D_p\lRa{\phi}\C P^1$ given by $u\mapsto \PP(T_pD_u)$.
This is just the Gauss map. Then
$S_1(Q)=\phi^{-1}[t]$. So, if $-p+l\in\phi^{-1}[t]$ then $S_1(Q)=\{-p+l,-p-l\}$. 
By Hurwitz (\cite[Cor.~IV.2.4]{Hart}) $\phi$
has a ramification divisor $R$ of degree 6. In particular, as $D_L$ is
irreducible, the reducible divisors in $|2\ell|$ which have a
tacnode are precisely $D_{-l}+D_l$, where $l\in D_L\setminus\T_2$. A similar
argument can be found in \cite[Lemma 5.1]{MukFF}.

\medskip
Let us turn now to $S_2(Q)$. 
\begin{prop}\label{p:Ri2Q}
For all $Q\in\Hilb^2\T$, with $Q=\{p,q\}$ (possibly $p=q$),
$R^0_2(Q)$ is a rank 2 vector bundle and
$R^1_2(Q)$ is a torsion sheaf. Then 
$R^1_2(Q)=\OO_{-p-q}$ and
$R^0_2(Q)=H_{S_1(Q)}\hat L^{-1}$.
In particular, $R^0_2(Q)$ is $\mu$-semistable.
\end{prop}
\begin{proof}
Consider $\{p\}\subset Q$ and use \seqref{e:RFofPstr} with $i=2$. In this case,
$R^1_2(P)=0$ and so if $R^1_2(Q)$
had non-zero rank it would be isomorphic to $\P_q$ which is impossible. Hence
$R^1_2(Q)$ has rank 0 and $R^0_2(Q)$ has rank 2 since $\chi(L^2\I_Q)=2$. 
This proves the first statement of the proposition.

Let $c_1(R^1_2(Q))=b$ and factor the middle map
of \seqref{e:RFofPstr} via $B\I_X$, where $X$ is a 0-scheme and $B$ is a line
bundle with $c_1(B)=-b$. Then $c_1(R^0_2(Q))=-2\ell+b$. But
\propref{p:R0FLI} implies that 
\[\mu(R^0_2(Q))=-2+b\cdot\ell/2<-1.\]
In other
words, $b\cdot\ell<2$. On the other hand, $b\cdot\ell\geq0$ from the
definition of $b$. Notice also that $B\I_X\subset \P_q$ so that
\[\irf2(B)\cong\irf2(B\I_X)\neq0.\]
This implies that $b^2\geq0$. Then there are two cases to consider:
\par\vspace{5pt}
\noindent\textbf{(1) $b\cdot\ell=1$.}\ \ Then the Hodge Index Theorem implies that
$b^2<1/2$
and so $b^2=0$. Hence, $B^*=\hat L_k\P_x$, say and the torus is a
product which we do not allow.
\par\vspace{5pt}
\noindent\textbf{(2) $b\cdot\ell=0$.}\ \ Since $b^2\geq0$ we must have $b^2=0$.
Hence, $B$ is flat. So $R^1_2(Q)=\OO_X$. Since
$H^0(L^2\I_Q)\neq0$ we can pick $L\hookrightarrow L^2\I_Q$. Let the quotient be
$A$. Then $A$ is supported on $D_L$ and
\[R^1_2(Q)\cong\rf1(A)\text{ and }
\rf0(A)=R^0_2(Q)/\hat L^*.\]
If $T\subset A$ is the torsion subsheaf of $A$
then $T$ is supported in codimension 0. On the other hand, \lemref{l:SubSch} implies that $T=0$.
Hence, $A=\OO_DL_\xhat$, for
some $\xhat$. The short exact sequence $\P_\xhat\to L_\xhat\to A$ implies that
\[\rf0(A)=\hat L^*\P_{-\xhat}\text{ and }\rf1(A)=\OO_{-\xhat}(=\OO_X).\]
Now apply
det to \seqref{e:RFofstr} to get $\P_\xhat=\det H_Q$ so that $\xhat=p+q$.
This deals with the two cases.

Observe that $\ch\bigl(R^0_2(Q)\tensor\hat L\bigr)=
(2,0,0)$ and since
$R^0_2(Q)$ is locally-free we must have $R^0_2(Q)\tensor \hat L=H_{\tilde Q}$
for some $\tilde Q\in\Hilb^2\dT$. The Mukai spectral sequence now
gives us a surjection $LH_{-\tilde Q}\to L^2\I_Q$ and so we see that
$\tilde Q=S_1(Q)$.
\end{proof}

\section{$|X|=3$}{}
We have dealt with $S_1(Y)$ in the preceding section.
We shall first treat the collinear case.
\begin{prop}\label{p:R12Ycol}
If $Y\subset D_v$ then $S_2(Y)=D_{v-\sum Y}$, $R^1_2(Y)$ is a degree
1, rank 1 torsion-free sheaf over $S_2(Y)$ and $R^0_2(Y)=\hat L^*\P_{-v}$.
Conversely, if $R^1_2(Y)$ is supported on a translate of
$D_L$ and torsion-free on its support then $Y$ is collinear.
\end{prop}
\begin{proof}
If $Y\subset D_v$ then we have a short exact sequence $L_{v}\to
L^2\I_Y\to A$, where $A$ is supported on $D_v$.  Apply $\muk$ to obtain
\[0\to \hat L^*\P_{-v}\to R^0_2(Y)\to\rf0(A)\to0\qquad\hbox{and}\qquad
R^1_2(Y)\cong\rf1(A).\]
We know from \lemref{l:SubSch} that the torsion of $A$ is supported in
dimension 1
and so $A$ is torsion-free over $D_v$. Since $\chi(A)=0$, Riemann-Roch implies
that $A$ has degree 1 over $D_v$.

Conversely, if $R^1_2(Y)$ is supported on $D_u$ then $R^0_1(Y)\cong\hat
L^*\P_a$, say and so the degree of $R^1_2(Y)$ is 1. On the other hand,
\seqref{e:MSSLI} shows that $\irf0(R^1_2(Y))=0$ and $Y$ is collinear.
\end{proof}
\begin{cor}\label{c:WITX}
If $X\in\Hilb^n\T$ with $n>3$ then $R^0_2(X)$ has rank 1 if and only if
$X$ is collinear. Hence, $L^2\I_X$ satisfies WIT$_1$ if and only
if $X$ is not collinear.
\end{cor}
\begin{proof}
Suppose that $X$ is collinear.
We must have $\rk(R^0_2(X))<2$ by \propref{p:Ri2Q}. On the other hand, if
$X\subset D_u$ then $X\subset D_u+D_v$ for any $v$ and so
$\Phi_2(X)_x\neq\varnothing$ for all $x$. Hence, $R^0_2(X)$ is a line bundle.
Conversely, if $R^0_2(X)$ is a line bundle then we induct on $n$. If
$X'\subset X$ has length $n-1$ and then \seqref{e:RFofPstr}
implies that
$R^0_2(X)\cong R^0_2(X')$ and $\irf2(R^0_2(X))\to L^2\I_{-X}$ implies that $X$
is collinear (by induction). This holds also if $n=4$ and so the induction
starts. 
\end{proof}
\begin{thm}\label{th:S2YeY}
Suppose that $Y\in Hilb^3\T$ is not collinear. Then $S_2(Y)\cong Y$
as schemes and $R^0_2(Y)\cong\hat L^{-2}\P_{-\Sigma Y}$.
\end{thm}
\begin{proof}
Let $Q\subset Y$ be a length 2 subscheme of $Y$. Then we can analyze
\seqref{e:RFofPstr}:
\[
0\lra R^0_2(Y)\lRa\phi R^0_2(Q)\lra\P_y
\lra R^1_2(Y)\lRa\psi\OO_{-p-q}\lra0.
\]
Let $A=\coker(\phi)$ and $B=\ker(\psi)$. We know from \propref{p:R0FLI} and
the \seqref{e:RFofPstr} for $Q\subset Y$ that
$R^0_2(Y)$ is a line bundle. Let $r=c_1(R^0_2(Y))$. Then
$\ch(A)=(1,-2\ell-r,2-r^2/2)$. Since $A$ is torsion-free it takes the form
$L^{-2}R^{-2}\I_X$ modulo $\Pic^0\T$ for some 0-scheme $X$. Then
$2-r^2/2=4+2r\cdot\ell+r^2/2-|X|$ and so $|X|=2+2r\cdot\ell+r^2\geq0$. From
the $\mu$-semistability of $R^0_2(Q)$ (see
\propref{p:Ri2Q}) the fact that and $\mu(R^0_2(Q))=-2$ we must have
$r\cdot\ell<-1$.
On the other hand, $c_1(B)\cdot \ell\geq0$. Hence,
$r\cdot\ell\geq-4$. We now treat the values of $r\cdot\ell$ separately.
\smallskip\noindent
\textbf{(1) $r\cdot\ell=-2$.}\ \ The Hodge Index Theorem implies that
$r^2\leq2$ but $-2+r^2=|X|\geq0$ and so $r^2=2$ with $|X|=0$. This implies that
$-r$ is a principal polarization and, since $-r\cdot\ell=2$, $r=\ell$. Then
$B$ is locally-free over a translate of $D_L$. But then $c_1(A)=\ell$
and so $c_1(R^1_2(Y))=\ell$. But, since $r(R^1_2(Y))=0$ we have that
$R^1_2(Y)$ is supported on a translate of $D_L$ and so by
\thmref{p:R12Ycol}, $Y$ must be collinear, a contradiction.

\smallskip\noindent
\textbf{(2) $r\cdot\ell=-4$.}\ \ In this case, $r^2\geq6$ and the Hodge
Index Theorem implies that $r^2=6$. Then the degree of $A$ is 0 and so, since
$A\to\P_y$ is non-zero, $A^{**}$ must be flat. Hence, $R=L^{-2}$ and
$X\in\Hilb^2\dT$. Then $B=\OO_X$ and
\[R^1_2(Y)=\OO_{S_2(Y)}\quad\text{with}\quad S_2(Y)\in\Hilb^3\dT.\]
If we apply det
to \seqref{e:RFofstr} then we see that $R^0_2(Y)=L^{-2}\P_{-p-q-y}$. If $Y$
is reduced then by repeating the argument above with each
$Q\subset Y$ we obtain
\[S_2(Y)=\{-p-q,-q-y,-y-p\}\]
and so is isomorphic to
$Y$. To discuss the others observe that by applying det to \seqref{e:RFofstr}
we see that
\[\sum S_2(Y)=-2\sum Y=-2p-2q-2y.\]
This deals with the case where the support of $Y$ contains at least 2 points.
Observe that the
argument used to eliminate case (2) also implies that
$|S_2(Y)|\subset\bigcup_{Q\subset Y}|S_2(Q)|$. If $Y$ is of type c, d or e
then the union of the $|S_2(Q)|$'s is a single point and therefore
$S_2(Y)$ is of type
c, d or e. Using the continuity
of
\[(p,q,y)\mapsto(-p-q,-q-y,-y-p)\] and $Y\to S_2(Y)$ we see that 
$Y\cong S_2(Y)$. In fact, we can use the sequences 
\begin{align*}
\{(0,0),\ (-x_n,0),\ (x_n,0)\}&\lra Y_c\\
\{(0,0),\ (x_n,0),\ (x_n,x_n^2)\}&\lra Y_d\\
\text{and}\quad\{(0,0),\ (0,x_n),\ (x_n,0)\}&\lra Y_e
\end{align*}
as $x_n\to0$.
Since the torus is an abelian Lie group
the addition of points corresponds to addition of these tangent vectors. Then
we see that $S_2(Y)$ for each of these is in the same configuration as
$x_n\to0$.
\end{proof}
\begin{note}
\begin{itemize}\item[]
\item[(1)] This theorem is the geometrical reason why there are no rank 2
stable bundles with Chern character $(2,0,-1)$. Any such bundle $E$
must fit into a short exact sequence $L_\xhat^{-1}\to E\to
L_{\xhat+\dhat}\tensor\I_Y$ for some $Y\in \Hilb^3\T$. But only collinear $Y$'s
give rise to locally-free extensions and for these
$H^0(E\tensor\P_\zhat)\neq0$ when $Y\subset \tau_{\xhat+\dhat+\zhat}D_L$. This
contradicts the stability of $E$.
\item[(2)] The case when $Y$ is collinear provides us with a good
example of a
sheaf where we know the Fourier transform cohomology sheaves $R^i $ and
even the boundary maps in the Mukai spectral sequence but we cannot
reconstruct the original sheaf from this alone. Whatever $Y$
is, the spectral sequence degenerates at the $E_2$ level. When $Y$ is
not collinear then $R^1_2(Y)=\OO_{\tilde Y}$ and
$R^0_2(Y)=L^{-2}\P_{-\tau}$ and from 
$\tau$ together with $\tilde Y$ we can easily reconstruct $Y$. But when $Y$ is
collinear $R^1$ is a degree 1 line bundle over $D_u$ and
$R^0_2(Y)=L^{-1}\P_{-\tau-u}$. Then $R^1$ is determined by its cohomology
jumping divisor in $\T$ which must be $D_v$. But applying $H^*$ to the long
exact sequence $R^0_2(Y)\to \Lshat\to\hat H_Y\to R^1_2(Y)$ we see that
$Y\subset D_v$
and so $v=\tau+u$. Hence, we only have the parameters $u$ and $\tau$ free
which are not enough to determine $Y$.
\item[(3)] There are precisely six 0-schemes $Y$ which are fat points of type c (so
  $Y\cong\spec\C[\epsilon]/(\epsilon^3)$) supported at $e\in\T$. These
  points are the ramification divisor of the Gauss map $\phi$
  discussed after Corollary \ref{c:genS1}. To see this observe that
  $R^1_2(Y_c)$ is properly torsion-free precisely when $L\I_{S_1(Q)}$
  satisfies the Cayley Bacharach condition, where $Q\subset Y_c$ is
  the unique length $2$ subscheme. But this happens precisely when when $S_1(Q)$
  is fat. So $\supp S_1(Q)\in D_L\cap\T_2$. This tell us that the
  inflection points of $D_L$ are the points of $D_L\cap\T_2$.
\end{itemize}\end{note}

\section{$|X|=4$}{}
From \corref{c:WITX} we can deduce the
following.
\begin{prop}
$L^2\I_Z$ satisfies WIT$_1$ if and only if $Z$ is not collinear. If $Z$
is collinear then the rank of $R^0_2(Z)$ is 1.
\end{prop}
\begin{prop}
If $Z\subset D_v$ then $R^0_2(Z)=\hat L^{-1}\P_{-v}$ and
$R^1_2(Z)=\hat L\P_{\sigma-v}\I_{2v-\sigma}$, where
$\sigma=\sum Z$.
\end{prop}
\begin{proof}
Suppose first that $Z$ contains some $Y$. Then
\seqref{e:RFofPstr} implies that $R^0_2(Z)\cong R^0_2(Y)=\hat L^{-1}\P_{-v}$.
We also have a
sequence $0\to\P_z\to R^1_2(Z)\to\OO_{D_{v-p-q-y}}(1)\to 0$, where
$(1)$ denotes
degree 1. This sequence cannot split so any torsion in $R^1_2(Z)$ is
supported on a proper dimension 1 subset of $D_{v-p-q-y}$. Sequence
\ref{e:RFofPstr} shows that $R^1_2(Z)$ satisfies WIT$_1$. If
$Y\subset Z$ then $R^1_2(Z)=\hat
L\P_x\I_a$. To compute $x$ we apply det to \seqref{e:RFofstr}.
This gives $x=\sigma-v$. To
compute $a$ we apply $\dmuk$ to give $L_{-v}\to L^2\I_{-Z}\to\OO_{D_{-x-a}}\P_a$
and so $-x-a=-v$, so $a=2v-\sigma$.

\end{proof}

\begin{thm}\label{t:R12ZisT}
If $Z\in\Hilb^4\T$ is not collinear then $T=R^1_2(Z)$ is a degree 3
torsion-free sheaf of rank 1 of a divisor $D'$ in
$|\hat L^2\P_\sigma|$, where $\sigma=\sum Z$. Moreover, $D'$ is a Kummer
divisor if and only if $Z$ contains a collinear length 3
subscheme. The restriction type of $D'$ is $(2,1)$ precisely when
there is exactly one such collinear length 3 subscheme and is $(1,1)$
if there are two.
\end{thm}
\begin{proof}
We know that if $Z$ is not collinear then $L^2\I_Z$ satisfies WIT$_1$.
$T$ has no torsion supported in dimension
0 since $T$ satisfies WIT$_1$.
If $Y\subset Z$ is not collinear then \seqref{e:RFofPstr} gives 
\[0\lra\hat L^{-2}\P_{-\tau}\lra\P_z\lra T\lra\OO_{S_2(Y)}\lra0.\]
The middle map factors through a degree 0 line bundle over a divisor in $|\hat
L^2\P_\sigma|$. 

Suppose that $Y\subset Z$ satisfies $Y\subset D_v$. Then \seqref{e:RFofPstr}
gives
\[0\lra\hat L^*\P_{-v}\lra\P_z\lra T\lra\OO_{D_{v-p-q-y}}(1)\lra0.\]
The middle map factors through $\OO_{D_{-v-z}}\P_z$ and so, $T$ is
supported on a Kummer divisor in $|\hat L^2\P_{v+z-v+p+q+y}|
=|\hat L^2\P_\sigma|$. Conversely, suppose that $T$ is supported on a Kummer
divisor $D_\alpha+D_\beta$. Notice that it is not possible for
$T$ to be reducible since
$\rf1(T)=L^2\I_{-Z}$. We need to understand the possible restriction types of
$T$ to $D_i$.
Since $T$ does not satisfy IT$_1$ there must be some
$X\in\Hilb^5\dT$ and $x\in\T$ such that $T=L^2\P_x\I_X/\OO$. Let $T_i$
be the restriction of $T$ to $D_i$ with degree $n_i$ and let
$T'_i=\ker(T\to T_i)$. Then $\chi(T'_i)=-n_i$. Since $\rf0(T)=0$ and $\rf1(T)$
is torsion-free we must have $\chi(T'_i)<0$ and so $n_i\geq1$. If $T$ is
locally-free over $D_1+D_2$ then $n_{2-i}=3-n_i$ and so $n_i\leq2$.
\propref{p:R12Ycol} implies that, if $T$ is not locally-free then
the restriction type $(n_1,n_2)=(1,1)$. This implies that the only possible
restriction types are $(2,1)$, $(1,2)$ and $(1,1)$. The first two cases
correspond to a single collinear $Y\subset Z$ and the latter to two different
collinear $Y\subset Z$.
\end{proof}

\begin{remark}\label{r:detofZ}
We can also see part of \thmref{t:R12ZisT} by observing that
$\chi(L^2\I_Z)=0$ and so $S_2(Z)$ is supported on a
divisor from $\widehat{L^2}$. We can be more precise about this by applying
$\muk$ to the structure sequence of $Z$. This implies that
$\det\muk(L^2\I_Z)^*\cong \Lhat^2\P_\sigma$.
\end{remark}

\section{$|X|=5$}{}

The case $|X|=5$ is a curious one because $R^1_2(X)$ can have torsion
reflecting the structure of $X$.

\begin{remark}
Since $D\cdot D_x=4$ for $D\in|2\ell|$ we see that if $X$ is collinear then $S(X)=\varnothing$.
\end{remark}
\begin{thm}\label{t:TorinW}
Suppose that $W\in\Hilb^5\T$ is not collinear. Let $S$ be the torsion subsheaf
of $R^1_2(W)$ and let $F=R^1_2(W)/S$. Then $S$ satisfies one of the following.

\begin{enumerate}
\item If $S=0$ then $R^1_2(W)=\widehat{L^2}\P_a\I_{W'}$ for some
$W'\in\Hilb^5\dT$. If $S$ is not zero then it must be a torsion-free sheaf of
rank 1 supported on a divisor in $|\hat L\P_x|$ for some $x$.
\item If $S$ is supported on a translate
of $D_L$ then it must have degree 0 or $-1$.
\item $S$ is a degree 0 line bundle on a translate
of $D_L$ if and only if $W$ contains a length 4 collinear subscheme.
\item $S$ is a degree $-1$ line bundle on a
translate of $D_L$ if and only if every length 4 subscheme
$Z\subset W$ contains a unique collinear length three subscheme but is
not itself collinear.
\end{enumerate}
\end{thm}
\begin{proof}
Recall from \corref{c:WITX} that $L^2\I_W$ satisfies WIT$_1$.
If $R^1_2(W)$ is torsion-free then from $\ch(R^1_2(W))=(1,2L,-1)$ we must have
$R^1_2(W)\cong\hat L^2\P_a\I_{W'}$ for some $W'\in\Hilb^5\dT$.
If $\OO_x\subset R^1_2(W)$ then this contradicts the fact that
$\irf0(R^1_2(W))=0$. This implies that the torsion of $R^1_2(W)$ is supported
on a divisor and is torsion-free over that divisor. This proves part
(1).

If we apply \seqref{e:RFofPstr} to
$Z\subset X$ then we obtain
\[0\lra R^0_2(Z)\lra\P_w\lra R^1_2(W)\lra R^1_2(Z)\lra0.\]
If $Z$ is collinear then $R^0_2(Z)\cong \hat L^*\P_v$ and $R^1_2(Z)$ is
torsion-free. This implies that $S=\P_w/L^*\P_v=\OO_{D_{v-w}}\P_w$.
If $Z$ is not collinear we have a commuting diagram with exact rows
and columns:
\begin{equation}\label{e:SinW}
\vcenter{\xymatrix{&&0\ar[d]&0\ar[d]&\\
&&S\ar@{=}[r]\ar[d]&S\ar[d]\\
0\ar[r]&\P_w\ar[r]\ar@{=}[d]&R^1_2(W)\ar[r]\ar[d]&R_2^1(Z)\ar[r]\ar[d]&0\\
0\ar[r]&\P_w\ar[r]&A\ar[r]\ar[d]&B\ar[r]\ar[d]&0\\
&&0&0&}}
\end{equation}
The bottom sequence implies that $B$ is torsion-free supported on a divisor.
In the case we are looking at, we have $R^1_2(Z)=T$, a torsion-free sheaf of
rank 1 supported on $D\in|L^2\P_a|$ for some $a$. We cannot have $S=T$ since
that would imply that $A=\P_w$ contradicting the fact that $R^1_2(W)$
satisfies WIT$_1$. Since $B$ is torsion-free on a divisor we must have that
$S$ is supported on $D_x\subset D$. This implies that $D$ is reducible and so
$Z$ contains a collinear $Y$. Moreover, the restriction type of $T$ is $(2,1)$
or $(1,1)$ and so the degree of $S$ is either 0 or $-1$. This proves (2) and
part of (3).

To complete (3) suppose that $S$ is degree 0 on $D_x$. Then
$\hat S\cong L\P_y\I_c$ for some
$x$ and $c$ and $R^1_2(W)/S\cong \hat L\P_a\I_b$. So if we apply $\invmuk$ to
\[0\lra S\lra R^0_2(W)\lra \hat L\P_a\I_b\lra0\]
we see that $L\P_y\I_c\subset L^2\I_{-W}$. This implies that
$W\setminus\{p\}\in\Hilb^4\T$ is collinear. 

We turn now to (4). If the degree of $S$ is $-1$ then no $Z\subset W$ is
collinear and the
restriction type of $R^1_2(Z)$ is $(2,1)$. So each length 4 subscheme $Z\subset W$ has a
unique collinear length 3 subscheme by Theorem
\ref{t:R12ZisT}. Conversely, suppose $L^2\I_W$ it WIT$_1$ and
its transform $R^1_2(W)$ is contained in a short exact sequence
\[0\lra S\lra R^1_2(W)\lra\hat L_x\lra0.\]
Note that $\dim\Ext^1(\hat L_x,S)=4$ and so there is a 9 dimensional family of
such extensions $R$. These extensions all satisfy WIT$_1$ and the
Fourier transform is generically torsion-free as both $\hat S$ and
$\irf0(\hat L_x)$ are locally-free, and any torsion in $\dmuk(R)$ must
be a quotient $\P_{\hat y}/L_x^{-1}$ coming from a lift of
$L_x^{-1}\to\hat S$. But the generic such map does not lift.
\end{proof}
\begin{note}\item[]
\begin{itemize}
\item[(1)] It is easy to see that if $R^1_2(W)$ is torsion-free then the factor of
$\P_x$ in it is given by $\P_{\Sigma W}$. 
\item[(2)] One might conjecture that $S_2(W)\cong W$ as schemes when $R^1_2(W)$ is
torsion-free. In fact, this need not be the case. As an example consider a
generic divisor $D\in|2\ell|$ and $D_x+D_{-x}\in|2\ell|$ also generic. Then we can
arrange that $D$ intersects $D_x$ and $D_{-x}$ transversely and that $D\cap
D_x\cap D_{-x}=\varnothing$. Now pick $Y\subset D\cap D_x$ and $Q\subset D\cap
D_{-x}$ and let $W=Y\cup Q$. Note that $W$ is not collinear because $Y$ is
collinear and $Q\not\subset D_x$. Then the line spanned by $D$ and
$D_x+D_{-x}$ in $|2\ell|$ is contained in $\Phi(W)_0$. Then the dimension of the
fibre of $R^1_2(W)\cong \hat L\P_x\I_{W'}$ at 0 is 3. This implies that $W'$
is not reduced at 0 and so cannot be isomorphic to $W$ which is simple by
construction. 
\end{itemize}\end{note}

\section{General $X$}{}
For $|X|>5$ we can prove that when $X$ is not
collinear then it has a torsion-free Fourier-Mukai transform if and only if it
does not contain a collinear colength 1 subscheme. The hard work goes into
proving the $|X|=6$ case and then the general case follows by induction.
\begin{thm}\label{t:R12XTF}
Let
$X\subset\T$ be a 0-dimensional subscheme of an irreducible ppas
$(\T,\ell)$ of length $|X|\geq6$. 
Suppose that $X$ is not collinear. Then $R^1_2(X)$ admits torsion $S$ if and
only if $X$ contains $X'$, a collinear length $|X|-1$ subscheme. Moreover, in
that case $S$ must be a degree 0 line bundle over some translate of $D_L$. In
particular, for a non-collinear length 5 subscheme $W\subset X$ we have
that the torsion subsheaf of $R^1_2(W)$ has degree 0.
\end{thm}
\begin{proof}
Observe that $L^2\I_X$ satisfies WIT$_1$ as $X$ is not collinear. This follows
from \corref{c:WITX}. If $X$ contains a collinear colength 1 subscheme $X'$
then we have an inclusion $L_\xhat\I_p\subset L^2\I_X$. This implies that
$R^1_1(\{p\})\subset R^1_2(X)$ and so $R^1_2(X)$ admits torsion.

Conversely, suppose that $R^1_2(X)$ contains torsion $S$. Arguing as in the
first paragraph of the proof of \thmref{t:TorinW}, we can deduce that $S$ must
be supported on a divisor and be torsion-free with rank one over that divisor.
We now induct on $|X|$. Consider first the case $|X|=6$. Suppose that we have
a non-collinear length 5 subscheme $W\subset X$. Then from \seqref{e:RFofPstr}
we see that $S\subset R^1_2(W)$. Then \thmref{t:TorinW} implies that
$\chi(S)\leq -1$. Let $F=R^1_2(X)/S$. Then
\[R^1_2(W)/S=F/\P_x,\]
where $\P_x=\ker(R^1_2(X)\to R^1_2(W))$. Let $G$ be the torsion subsheaf of
$R^1_2(W)/S$ (regarded as a sheaf on its support so that $G$ is supported at
points). Then $G\subset F/\P_x$. But $\Ext^1(G,\P_x)=0$ and so this map lifts
to a map $G\to F$. But $F$ is torsion-free (as an $\OO_\T$-module) and this
implies that $G=0$. In other words, $S$ is the torsion subsheaf of $R^1_2(W)$
as well. Then \thmref{t:TorinW} implies that $\chi(S)=-1$ or $-2$.

Suppose that $\chi(S)=-2$. Then $F$ is a rank two vector bundle and
$F/\P_x=R^1_2(W)/S=\hat L_y$ for some $y$. This implies that
$\irf2(F)=\OO_{-x}$ --- apply $\invmuk$ to the short exact sequence
\[0\lra\P_x\lra F\lra\hat L_y\lra0.\]
But this contradicts the fact that $R^1_2(W)$ satisfies WIT$_1$.

Then $\chi(S)=-1$ and we have a short exact sequence
\[0\lra\P_x\lra F\lra\hat L_y\I_p\lra0.\]
This implies that $\irf0(F)=0$ and so $F$ satisfies WIT$_1$. We then have a
short exact sequence
\[0\lra\hat S\lra L^2\I_{-X}\lra\hat F\lra0.\]
But $\hat S\cong L_{\hat x}\I_p$ for some $\xhat$ and $p$ and so
$X\setminus\{-p\}$ is collinear.
This proves the theorem for the case $|X|=6$.

Now suppose that $|X|>6$ and we have proved the theorem for all 0-schemes of
length $|X|-1\geq6$. We suppose that $R^1_2(X)$ contains torsion $S$ and
suppose that a length $|X|-1$ subscheme $X'\subset X$ is not collinear. Then
$S\subset R^1_2(X')$ and the induction hypothesis implies that there is a
length $|X'|-1$ subscheme $X''\subset X'$ which is collinear. Let
$Q=X\setminus X''$. then we have a short exact sequence
\begin{equation}
0\lra L_\xhat\I_Q\lra L^2\I_X\lra A\lra0
\label{e:QinX}
\end{equation}
for some $\xhat$ and line bundle $A$ supported on $D_\xhat$. Note that
$\chi(A)=5-|X|$. 

Suppose, for a contradiction, that $A$ is locally-free on its support. Then,
since $\deg(A)<0$ we must have that $A$ satisfies IT$_1$. Then applying $\muk$
to \ref{e:QinX} we have a short exact sequence
\[0\lra\hat L_x\I_{Q'}\lra R^1_2(X)\lra\hat A\lra0,\]
for some $x$. But $R^1_2(X)$ admits torsion. This is a contradiction.

So $\OO_{\tilde X}\subset A$. Let $\tilde A=A/\OO_{\tilde X}$ be locally-free on
its support. Then $\chi(\tilde A)<\chi(A)$. But there is a surjection
$L^2\I_X\to \tilde A$ and its kernel is torsion-free of rank 1 with
singularity set strictly contained in $Q$. This implies that it is the ideal
sheaf of a single point $\{p\}$ and we have $L_x\I_p\subset L^2\I_X$.

We now see that $S$ is either 0 or has degree 0 over a translate of $D_L$.
This means that if $W\subset X$ is non-collinear then $S\subset R^1_2(W)$ and
so the degree of the torsion subsheaf of $R^1_2(W)$ must be greater or equal
to 0.
\end{proof}

We know that $D\cdot D'=8$ for $D,D'\in|2\ell|$. In particular, this means that
there must exist $X\in\Hilb^8\T$ such that $S_2(X)\in\Hilb^2\dT$. We know from
\thmref{t:TorinW} that when $R^1_2(W)$ is torsion-free then
$S_2(W)\in\Hilb^5\dT$. We would now like to compute the length of $S_2(X)$ for
$|X|\geq6$. This is answered in the following theorem.
\begin{thm}
Let
$X\in\Hilb^n\T$ with $n\geq6$. Suppose that $X$ is not collinear and that
$R^1_2(X)$ is torsion-free. Then 
\[\length\bigl(\sing(R^1_2(X))\bigr)=\length(S_2(X))\leq3.\]
Furthermore, if $|X|\geq7$ then $\length(S_2(X))\leq2$.
\label{t:LthofX}
\end{thm}

\begin{proof} Observe first that 
\begin{equation}\label{e:R2R12X}
\irf2(R^1_2(X)^{**})=0.
\end{equation}
This follows by applying $\invmuk$ to
\[0\lra R^1_2(X)\lra R^1_2(X)^{**}\lra\OO_{S_2(X)}\lra0.\]

Let $|X|=6$. Consider $W\subset X$. By \thmref{t:R12XTF}, we know that $W$
cannot be collinear. Suppose first that $R^1_2(W)$ is torsion-free. Without
loss of generality, assume that $R^1_2(W)=\hat L^2\I_{W'}$. Then we have the
following 3 by 3 diagram:
\begin{equation}\label{e:R12DD}
\vcenter{\xymatrix{&&0\ar[d]&0\ar[d]&\\
0\ar[r]&\P_v\ar[r]\ar@{=}[d]&R^1_2(X)\ar[r]\ar[d]&L^2\I_{W'}\ar[d]\ar[r]&0\\
0\ar[r]&\P_v\ar[r]&R^1_2(X)^{**}\ar[d]\ar[r]&L^2\I_{W'\setminus
  S_2(X)}\ar[d]\ar[r]&0\\
&&\OO_{S_2(X)}\ar@{=}[r]\ar[d]&\OO_{S_2(X)}\ar[d]\\
&&0&0&}}
\end{equation}
The middle horizontal sequence cannot split by \ref{e:R2R12X}. But
\[H^1(L^2\I_{W'\setminus S_2(X)})=0\quad\text{if}
\quad|W\setminus S_2(X)|\leq1.\]
So $5-\length(S_2(X))\geq2$.

Now suppose that $R^1_2(W)$ has non-zero torsion subsheaf $S$. By
\thmref{t:TorinW} and the last part of \thmref{t:R12XTF} we have
\[R^1_2(W)/S\cong\hat L_x\I_p,\]
for some $x$ and $p$. Then
\[\ker(R^1_2(X)\to\hat L_x\I_p)=\hat L_y\I_Q,\quad Q\in\Hilb^2\dT\]
and so we have $Q\subset S_2(X)\subset Q\cup\{p\}$ and we conclude that
\[2\leq\length(S_2(X))\leq3.\]

This completes the proof of the first part of the theorem.

Suppose now that $|X|=7$. We may assume as before that $X$ contains a
non-collinear length 6 subscheme $X'$. By \thmref{t:R12XTF} we have two cases
to consider: either $R^1_2(X')$ is torsion-free or it contains a torsion sheaf
$S$ with $E=R^1_2(X')/S$ a torsion-free sheaf with Chern character
$(2,\ell,0)$.

If $R^1_2(X')$ is torsion-free then we can argue as before. Let
$F=R^1_2(X)^{**}/\P_x$, where $\P_x=\ker\bigl(R^1_2(X)\to R^1_2(X')\bigr)$.
Then $\ch(F\tensor L^*)=(2,0,0)$. On the other hand, $F$ is locally-free and
so $F$ must be a homogeneous bundle. Therefore,
$\irf2(R^1_2(X)^{**})\neq0$. This implies that $\irf2(R^1_2(X))\neq0$;
a contradiction. 

On the other hand, suppose that $R^1_2(X')$ admits torsion.
We now need a lemma:
\begin{lemma}\label{l:R12SLF}
If $X'\in\Hilb^6\T$ is not-collinear and $R^1_2(X')$ admits torsion $S$ then
$E=R^1_2(X')/S$ is locally-free.
\end{lemma}
\begin{proof} 
As above, if $W\subset X'$ is not collinear then $R^1_2(W)/\mathrm{torsion}=\hat L_x\I_p$
and we have a short exact sequence
\begin{equation}\label{e:PtoE}
0\lra\P_x\lra E\lra \hat L_x\I_p\lra0.
\end{equation}
Suppose that $E$ is not locally-free then
we have the following 3 by 3 diagram:
\[\xymatrix{&&0\ar[d]&0\ar[d]&\\
0\ar[r]&\P_y\ar[r]\ar@{=}[d]&E\ar[d]\ar[r]&\hat
L_x\I_p\ar[r]\ar[d]&0\\
0\ar[r]&\P_y\ar[r]&E^{**}\ar[r]\ar[d]&L_x\ar[r]\ar[d]&0\\
&&\OO_p\ar@{=}[r]\ar[d]&\OO_p\ar[d]\\
&&0&0}\]
The middle horizontal sequence is split because $\hat L^{-1}$ satisfies IT$_2$ but
this implies that there is a non-zero map $E\to\P_x$ which contradicts the
fact that $\irf2(E)=0$. Alternatively, \seqref{e:PtoE} cannot be split
and so $\{p\}$ satisfied the Cayley-Bacharach condition with respect
to $|\hat L|$.

On the other hand, if all $W\subset X'$ are collinear then we have $E\cong
R^1_2(W)$ and so $E$ is locally-free by \corref{c:WITX}.
\end{proof}

Returning to the proof of the theorem, we now have a short exact sequence
\[0\lra \hat L_z\I_Q\lra R^1_2(X)\lra E\lra0\]
for some $z$ and $Q\in\Hilb^2\dT$ and locally-free sheaf $E$. This implies
that
\[\length\bigl(\sing(R^1_2(X))\bigr)=2.\]
This completes the proof of the theorem.
\end{proof}

\begin{cor}
If $R^1_2(X)$ is not torsion-free for some $X\in\Hilb^n\T$, $n\geq6$, then
$\Phi_2(X)_\xhat\subset Ks_\xhat(\T)$ for all $\xhat\in S_2(X)$.
\end{cor}
\begin{proof}
By \thmref{t:R12XTF} we know that $X$ contains a collinear length $n-1$
subscheme $X'\subset D_y$. Let $X=X'\cup\{x\}$. Then $X\subset D_y+D_z$ for
$z\in S_1(\{x\})$. 
Let the torsion subsheaf of $R^1_2(X)$ be $S$ as usual.
Note that $S_2(X)$ equals $D_\yhat$ as a scheme for some $\yhat$. This is
because $R^1_2(X)/S$ is locally-free. This means that
\[\Phi_2(X)=\bigl\{D_y+D_{-z}\mid z\in S_1(\{x\})\bigr\}\]
as required
\end{proof}

\section{Stability Properties of $R^1_i(X)$}
We now turn to stability properties of $R^1_i(X)$ for $i=1$ and $i=2$.
The following has been proved before (Mukai \cite[Thm 0.3]{MukFF}
and Yoshioka \cite[Prop 3.5]{Yosh01}). We give a direct and brief proof in the spirit of this paper.
\begin{thm}
If $|X|\geq2$ then $R^1_1(X)$ is $\mu$-stable.
\end{thm}
\begin{proof} We induct on the length $|X|$.
Observe that the $|X|=2$ case is trivial because $R^1_1(X)$ is rank $1$
torsion-free and so $\mu$-stable. Now assume that $|X|>2$. Pick $X'\subset
X$ of length $|X|-1$ and let $X=X'\cup \{x\}$ as schemes. Then we have a short exact sequence,
\[0\to \P_x\to R^1_1(X)\to R^1_1(X')\to 0.\]
Now that slope of $R^1_1(X)$ is $2/(|X|-1)$ and
$\mu(R^1_1(X'))=2/(|X|-2)$. Suppose $M$ is maximally
destabilizing. Then the map $M\to R^1_1(X)$ cannot lift to $\P_x$ and
so we have a non-trivial map $M\to R^1_1(X')$. If $r(M)<|X|-2$ then
the induction hypothesis implies that $d(M)<2$ which is impossible. On
the other hand, if $r(M)=|X|-2$, we have that the map $M\to R^1_1(X')$
must inject and so $\chi(M)\leq-1$. But
$\ch(R^1_1(X)/M)=(1,0,-1-\chi(M))$ and, by the maximality of the
destabilizing sheaf, $R^1_1(X)/M$ is torsion-free and so
$\chi(M)=-1$. But that is a contradiction as $M$ would then split the
original short exact sequence.
\end{proof}
Conversely, given a $\mu$-stable sheaf $E$ of Chern character
$(r,\ell,-1)$ with $r\geq 2$ we see that it cannot be IT$_1$ and it
must admit an injection from a homogeneous bundle of rank $r-1$. The
quotient must take the form $L_x\I_Q$ and so $E$ must be WIT$_1$ with
transform given by $L_y\I_X$ for some $0$-scheme $X$ of length
$r+1$. So we recover:
\begin{cor}\cite[(Thm 0.3 and Cor 5.5)]{MukFF}\label{c:muk2}
The moduli space $\mathcal{M}(r,\ell,-1)$ of $\mu$-stable sheaves is
compact and is biholomorphic
to $\Hilb^{r+1}\T\times\hat\T$ for $r\geq2$. The
non-locally free boundary has transform given by a collinear $0$-scheme.
\end{cor}

We now turn to the case $i=2$. The interesting cases are $|X|\geq6$.
Observe first that if $X$ is collinear then $R^0_2(X)=L^{-1}\P_y$ and
$R^1_2(X)$ has rank $|X|-3$ and degree $2$. It is an easy induction to
show that $R^1_2(X)$ must always be $\mu$-stable and locally-free.

We now look in detail at $|X|=6$. This is also treated, from the
point of view of moduli spaces, in \ moduli spaces, in \cite{Mac0} (see, in particular,
Theorem 4.10).
\begin{thm}\label{t:stab6}
 Let $X\in\Hilb^6\T$. Suppose that $X$
is not collinear.\begin{enumerate}
\item If $R^1_2(X)$ is locally-free then $R^1_2(X)$ is
$\mu$-stable.
\item If $R^1_2(X)$ is torsion-free then it is Gieseker stable.
\begin{enumerate}
\item $X$ contains a collinear length 4 subscheme if and only if
  $R^1_2(X)$ is $\mu$-destabilized by $L_x\I_Q$ for some $x$ and
  $Q$. Then $\sing(R^1_2(X)|=2$.
\item Otherwise $R^1_2(X)$ may be destabilized by $L_x\I_Y$, for some
  $x$ and $Y$. Then $|\sing(R^1_2(X)|=3$.
\end{enumerate}
\item If $R^1_2(X)$ admits torsion $S$ then $R^1_2(X)/S$ is
$\mu$-stable.
\end{enumerate}
\end{thm}
\begin{proof}
Suppose first that $R^1_2(X)$ is torsion-free. Then by \thmref{t:R12XTF} we
have a non-collinear length 5 subscheme $W\subset X$ and a short exact
sequence from \seqref{e:RFofPstr}
\[0\lra\P_x\lra R^1_2(X)\lra R^1_2(W)\lra0\]
Note that $\mu(R^1_2(X))=2$. Suppose that $M\subset R^1_2(X)$ is a
rank $1$ destabilizing subsheaf. Then $M$ takes the form
$L^n_y\I_{X'}$ for some $n>0$ and 0-scheme $X'$. Then the composite to
$R^1_2(W)=L^2_z\I_{W'}$ must inject. Hence, $\P_x\to R^1_2(X)/M$ also
injects and so $n\leq 2$. If $n=2$ then $R^1_2(X)/M\cong\P_x$ and so
$|X'|=3$. But $W\subseteq X'$ which is a contradiction. So we are left
with the case $n=1$. We may assume the quotient $R^1_2(X)/M$ is
torsion-free and so takes the form $L_w\I_{X''}$. The Euler characters
now tell us that $|X'|+|X''|=3$. If $|X'|=0$ then $M$ is IT$_0$ and
this contradicts the fact that $R^1_2(X)$ is WIT$_1$. Otherwise,
$R^1_2(X)$ must have a non-zero singularity set (containing
$X'$). This establishes the first part.

If $|X'|=1$ then the transform of $M\to R^1_2(X)$ gives a non-zero map
$\hat M\to L^2\I_X$. But this is impossible as $\hat M$ is a torsion
sheaf. Hence $|X'|>1$ and so $R^1_2(X)$ is Gieseker stable (as
$\chi(M)/1\leq -1$ while $\chi(R^1_2(X))=-1/2$).

Observe that $X$ has a collinear length $4$ subscheme if and only if
we have $L_y\I_{Q'}\to L^2\I_X$ for some $y$ and $Q'$. This happens if and
only if there is an injection $\widehat{L_y\I_{Q'}}\to R^1_2(X)$. But
$\widehat{L_y\I_{Q'}}=L_x\I_Q$ for some $x$ and $Q$.

Now assume $R^1_2(X)$ has torsion sheaf $T$ and let $E=R^1_2(X)/T$. We
have already seen that $E$ must be locally-free. Then we have a short exact sequence
\[0\to\P_y\to E\to L_x\I_p\to0\]
for points $x$, $y$ and $p$. But now any $\mu$-destabilising bundle would
have to have rank 1 and slope at least $2$. But then it cannot map to
either $L_x\I_p$ or to $\P_y$. This establishes the final part.
\end{proof}
%In particular, if $E\to R^1_2(X)$ is a map from a $\mu$-semistable
%sheaf $E$, we must have $d(E)\leq 2$.

\section{Applications}
We can deduce the following result (which can also be deduced from
classical projective geometry)
\begin{thm}
Suppose $D'$ is a genus five curve in an irreducible principally
polarized abelian variety $(\T,\ell)$.
If $D'$ is smooth then $D'$ has no $g^1_3$. There are (singular) irreducible
divisors $D'\in|2\ell|$ which admit a 1-dimensional family of $g^1_3$'s.
\end{thm}

\begin{proof}
Observe first that if $T$ is a degree 3 line bundle over a smooth divisor
$D'\in|2\ell|$ and if $T$ is the restriction of $L_\xhat\I_p$ then $T$
has no $g^1_3$. This is because the exact sequence
\[0\lra L^*\P_{\hat y}\lra L_\xhat\I_p\lra T\lra0\]
shows that $H^0(T)\cong H^0(L_\xhat\I_p)$ which has dimension 0 or 1.
This means that if $D'$ is not reducible and admits a $g^1_3$ then that
$g^1_3$ must satisfy WIT$_1$ with transform $L^2\P_\xhat\I_Z$ for some
$Z\in\Hilb^4\T$.

Observe that
\[\dim(H^0(T))=\dim\bigl((L^2\P_\xhat\I_Z)_{(0)}\bigr),\]
where
$A_{(x)}$ denotes the fibre of a sheaf $A$ at $x$. Then the
dimension of $H^0(T)$ is greater than 1 only if $(\I_Z)_0$ is not a regular
ideal. This can only happen if $Z$ contains either $Y_d$ or $Y_e$. But by
\thmref{th:S2YeY} we see that $S_2(Y)$ has type d or e. On the other hand,
$S_2(Y)\subset S_2(Z)$ and so $S_2(Z)$ must be singular. If we take the
example
\[Z\cong\C[\epsilon,\eta]/(\epsilon^3,\eta^2,\epsilon\eta)\]
then we see that $Z\subset S_2(Z)$ by \thmref{th:S2YeY}. The subset of
$\Hilb^4\T\cross\dT$ of such $Z$ is 6-dimensional, whereas the space of
singular divisors in $\PP\Lshat$ is at most 4-dimensional. This implies that
$\dim(W^1_3)\geq1$ for such $S_2(Z)$. On the other hand, such $S_2(Z)$ cannot
be generically reducible since then every $Y_c\subset\T$ would be collinear
which is impossible because it would tell us that $Y_c$ is, generically,
contained in two different translates of $D_L$. This contradicts
\propref{p:R12Ycol}. 
\end{proof}

We can also use our results to understand the set of singular divisors
in $|2\ell|$. Recall that $\T_2$ denotes the $2$-torsion points of $\T$.
\begin{thm}
 Let
$\Sigma_x$ denote the linear system of divisors in $|2\ell|$ which
pass through $x$. \begin{enumerate}
\item If $x\not\in\T_2$ then the only singular divisors in $\Sigma_x$ are
Kummer divisors.
\item If $x\in\T_2$ then any divisor $D$ in $\Sigma_x$ is
singular. If $D$ is not a Kummer divisor then its singularity set
consists of at most 3 nodes. There are $120$ (skew) lines of divisors with 2
nodes and each line contains fourteen points whose corresponding divisors have exactly 3 nodes.
\end{enumerate}
In particular, the locus of singular divisors consists of $Ks(\T)$
union all sixteen singular planes in $\PP^3=|2\ell|$.
\end{thm}
\begin{proof}
Observe that if $Q$ is a length 2 subscheme supported at $x$ then $S_2(Q)=\{-2x\}$.
(1) If $x\not\in\T_2$ then $h^0(L^2\I_Q)=1$. Let $Q_1$ and
$Q_2$ be two distinct length 2 subschemes supported at $x$ and suppose that
$\Phi_2(Q_1)_0\cap\Phi_2(Q_2)_0$ is 1 dimensional. Then this holds for
all pairs $Q_1$ and $Q_2$ since each divisor in the intersection has a
singularity at $x$. But then such singular divisors must also contain
$Y=Q_1\cup_x Q_2$ for some $Q_1\neq Q_2$. But $S_2(Y)=\{-2x\}$ and so
$\Phi_2(Y)_0$ consists of
a single divisor. This is a contradiction and shows that there is at
most 1 singular divisor in $\Sigma_x$ with a singularity at $x$. But
for any $x$ we can find $y$ such that $D_y+D_{-y}$ has a node at $x$.

(2) If $x\in\T_2$ then $\Phi_2(Q)_0=\Sigma_x$ and so every divisor in
$\Sigma_x$ is singular at $x$. Note that there is a line of Kummer
divisors in $\Sigma_x$ which have a tacnode or are multiple. The sets
$\Sigma_x$ correspond to singular planes in $\PP^3$. These intersect
along lines and any three of them intersect in a single point.
In fact the non-Kummer points correspond to G\"opel tetrahedra (see
\cite{BirkLange}). Since there are no degenerate G\"opel tetrahedra, we see
that these points are all distinct (so the divisors cannot have more
than 4 singular points unless they are multiple Kummer divisors.

It remains to prove that if a divisor $D\in\Sigma_x$ is not a Kummer
divisor then it has only nodes for singularities.
Observe now that $\Phi_2(Y)_0$ has
dimension 1 for all length 3 $Y$ supported at $x$. Pick a regular $Y$
and $D\in\Sigma_x\setminus\Phi_2(Y)_0$. Let $Y'$ be a length 3
subscheme supported at $x$ and contained in $D$. Then
$\Phi_2(Y)_0\cap\Phi_2(Y')_0=\{D\}$.
The space of pairs of regular length $3$ subschemes supported at $x$ is
given by $S^2\PP^1\cong\PP^2$ and the above observation furnishes us
with a map $\PP^2\setminus\Delta\to\Sigma_x\cong\PP^2$, where $\Delta$
denotes the diagonal. This must be
finite since otherwise there would be divisor containing all regular
length $3$ subschemes supported at $x$ and this is impossible. On the
other hand the map has degree one since it is given by the
intersection of hypersurfaces and so it must be the identity. The
extension over the diagonal maps to the Kummers and so we see that the
non-Kummer divisors have a simple node at $x$.
\end{proof}

\begin{remark} We can also see directly that there is at most one
divisor in $|2\ell|$ which is singular at three given distinct points
because for a generic choice of length $6$ subscheme $X$ supported at the
three points $h^1(L^2\I_X)\leq1$ because $X$ does not contain collinear
length $5$ subschemes and is not, itself, collinear.
\end{remark}

\section*{Appendix}{}
We will now summarize the computations of the Fourier-Mukai transforms derived
in section 5--10.

\begin{center}\begin{tabular}{@{}|c|l|l|l|l|l|l|@{}}
\hline
$X$\bigstrut&$S_1(X)$&$S_2(X)$&$R^1_1$&$R^0_2$&$R^1_2$&where\\
\hline
$P$\bigstrut&$D_p$&$\varnothing$&$\P_{p}\OO_{D_{-p}}$&(Note 1)&$0$&\\
\hline
$Q$\bigstrut&${\scriptstyle\{-p+l,-p-l'\}}$&$\{-p-q\}$&$\I_{S_1}\Lhat$&&$\OO_{S_2}$&$p-q=l-l'$\\
\hline
$Y$\bigstrut[t]&$\varnothing$&$\cong Y$&&$\Lhat^{-2}\P_{-\tau}$&$\OO_{S_2}$&$\tau=\sum
Y$\\
&\bigstrut $\{v\}$&$D_u$&&$\Lhat_{-v}^{-1}$&(Note 2)&$u=v-\tau$\\
$Z$\bigstrut &$\varnothing$&$D'$&&$0$&(Note 3)&$D'\in|\Lhat^2\P_{\sigma}|$\\
&$\{v\}$\bigstrut[b]&$\{2v-\sigma\}$&&$\hat L_{-v}^{-1}$&$\Lhat_{\sigma-v}\I_{2v-\sigma}$
&$\sigma=\sum Z$\\
\hline
\end{tabular}\end{center}

\begin{enumerate}
\item $R^1_2(P)$ is a rank $3$ $\mu$-stable vector bundle.
\item $R^1_2(Y)$ is a degree 1 line bundle over $D_u$.
\item $R^1_2(Z)$ is a degree $3$ line bundle over $D'$
\end{enumerate}

\begin{center}\begin{tabular}{@{}|c|l|l|l|l|l|@{}}
\hline
$X$\bigstrut &$S_1(X)$&$S_2(X)$&$R^0_2$&$R^1_2$&where\\
\hline
$W$\bigstrut[t]&$\varnothing$&$W'$&$0$&$\hat L^2\P_\alpha\I_{W'}$&$W'\in\Hilb^5\dT$\\
&\bigstrut $\varnothing$&$D_u$&$0$&$T\ltimes\Lhat_\xhat\I_y$&collinear
$Z\subset W$, $\deg(T)=0$ on $D_{\hat u}$\\
&$\varnothing$\bigstrut &$D_u$&$0$&$T\ltimes\Lhat_\xhat$&$\exists x$,
$\forall Z\subset W$, $\exists Y\subset Z$, $Y\subset D_x$\\
&&&&&$\deg(T)=-1$  on $D_{\hat u}$\\
&\bigstrut[b]$\{v\}$&$\varnothing$&$\Lhat^{-1}_v$&(Note 1)&\\
\hline
\bigstrut[t]$X$&$\varnothing$&$\varnothing$&$0$&(Note 2)&$|X|\geq6$\\
&\bigstrut $\varnothing$&$\{u\}$&$0$&(Note 3)&$|X|\geq6$\\
&\bigstrut $\varnothing$&$\{u,v\}$&$0$&(Note 3)&$6\leq|X|\leq8$\\
&\bigstrut $\varnothing$&$\{u,v,z\}$&$0$&(Note 3)&$|X|=6$\\
&\bigstrut $\varnothing$&$D_u$&$0$&$T\ltimes E $&$|X|\geq6$, collinear $X'\subset X$ \\
&&&&&\bigstrut colength $1$, $\deg(T)=0$ on $D_u$\\
&$\{v\}$&\bigstrut[b]$\varnothing$&$\Lhat^{-1}_v$&(Note 4)&$|X|>5$\\
\hline
\end{tabular}\end{center}

\begin{enumerate}
\item $R^1_2(W)$ is a rank 2 $\mu$-stable vector bundle over $\dT$ with Chern
character $(2,\ell,0)$.
\item $R^1_2(X)$ is a vector bundle over $\dT$ with Chern character
$(|X|-4,2\ell,-1)$.
\item $R^1_2(X)$ is a torsion-free sheaf over $\dT$ with Chern character
$(|X|-4,2\ell,-1)$ and singularity set equal to $S_2(X)$. This is
$\mu$-semistable when $|X|=6$.
\item $R^1_2(X)$ is a vector bundle over $\dT$ with Chern character
$(|X|-3,\ell,0)$.
\end{enumerate}

\bibliographystyle{amsplain}
\bibliography{p51}

\end{document}